\documentclass[twoside,12pt]{article}
\makeatletter
\def\ifundefined#1{\expandafter\ifx\csname#1\endcsname\relax}
\RequirePackage{theorem}

  \usepackage[breaklinks=true,%
  bookmarks=true,%
  backref=page,%
  pagebackref=true
]{hyperref}
  \AtBeginDocument{\hypersetup{breaklinks=true,%
      bookmarks=true,%
      backref=page,%
      pagebackref=true
    }}
    \AtBeginDocument{
   \makeatletter
   \pagestyle{myheadings}
   \markboth{\hfill\ifundefined{@authorshort} \@author
                   \else \@authorshort 
                   \fi\hfill}
            {\hfill\ifundefined{@titleshort} \@title 
                   \else \@titleshort
                   \fi\hfill}
   \makeatother
 }
        \providecommand{\MR}[1]{MR\href{http://www.ams.org/mathscinet-getitem?mr=#1}{#1}}
  \providecommand{\Zbl}[1]{Zbl\href{http://www.emis.de:80/cgi-bin/zmen/ZMATH/en/zmathf.html?first=1&maxdocs=3&type=html&an=#1&format=complete}{#1}}
     \theorembodyfont{\slshape}
        \newtheorem{thm}{Theorem}[section]

     \newtheorem{cor}[thm]{Corollary}

     \theorembodyfont{\upshape}
     \newtheorem{defn}[thm]{Definition}
     
     \newtheorem{example}[thm]{Example}
     \newtheorem{conj}[thm]{Conjecture}

     \newtheorem{rem}[thm]{\mdseries\scshape Remark}
\ifundefined{proof}

\newcommand{\proofname}{Proof}
\fi

\providecommand{\dedicatory}[1]{}
\providecommand{\keywords}[1]{\begingroup \def \protect {\noexpand \protect \noexpand }\xdef \@thefnmark { }\endgroup \@footnotetext{{\em Keywords and phrases.\/} #1}}
\providecommand{\AMSMSC}[2]{\begingroup \def \protect {\noexpand \protect \noexpand }\xdef \@thefnmark { }\endgroup \@footnotetext{{1991 \it 
Mathematical Subject Classification.\/} Primary: #1; Secondary: #2.}}

\newcommand{\authorshort}[1]{\gdef\@authorshort{#1}}
\newcommand{\titleshort}[1]{\gdef\@titleshort{#1}}
   {\@definecounter{equation}}{\@newctr {equation}[section]}
   
\def\p@enumi{\thethm.}

\def\p@enumi{}

\newcommand{\ppnum}[3]{\gdef\rep@rtenum{#2}
\gdef\rep@rteyear{#3}\gdef\wh@reappear{#1}}
   \let\@ldmaketitle=\maketitle
   \renewcommand{\maketitle}{ \vspace*{-2cm}  \makeatletter
   {\def\newpage{}
   {\footnotesize\noindent
   \parbox[t]{0.4\textwidth}{\noindent Preprint
   \texttt{\rep@rtenum}, \rep@rteyear\\
   \wh@reappear} \hfill
   \parbox[t]{0.5\textwidth}{\whereappear}}
   \@ldmaketitle}\makeatother}
   \newcommand{\whereappear}{} % Journal or Proceedings data
\ppnum{}{clf-alg/kisi9602}{1996}
\renewcommand{\whereappear}{Printed in\\
V.~Dietrich, K.~Habetha, and G.~Jank (eds): 
\emph{Clifford algebras and their application in mathematical physics. 
Aachen 1996}.
Kluwer Academic Publishers, 1998, pp. 175-184.}
\newcommand{\comment}[1]{}

\usepackage{amsfonts}

\newcommand{\algebra}[1]{\ensuremath{\mathfrak{#1}}}

\newcommand{\object}[2][\,]{\ensuremath{\mathrm{#2}#1}}
\newcommand{\Space}[2]{\ensuremath{ {\mathbb{#1}^{#2}} }}

\newcommand{\FSpace}[2]{{\ensuremath{ #1_{#2} }}}

\ifundefined{qed}
    \DeclareMathSymbol{\qed}{0}{AMSa}{"03}
\fi

\newcommand{\modulus}[1]{\left|#1\right|}
\newcommand{\scalar}[2]{\langle #1,#2\rangle}

\providecommand{\eqref}[1]{\textup{(\ref{#1})}}

\newcommand{\person}[1]{\textsc{#1}}

\newcommand{\matr}[4]{{\ensuremath{ \left( \begin{array}{cc}
#1 & #2 \\ #3 & #4
\end{array}\right) }}}
\makeatother
\gdef\p@enumi{}

\titleshort{How Many Essentially Different Functional Theories Exist?}
\authorshort{Vladimir V. Kisil}
 
\begin{document}
\title{How Many Essentially Different\\ Functional Theories Exist?}
\author{Vladimir V. Kisil\\
               Institute of Mathematics,\\
               Economics and Mechanics,
               Odessa State University\\
               ul. Petra Velikogo, 2, 
               Odessa}
\date{June 12, 1996}
\maketitle
\begin{abstract}
The question in the title is ambiguous. At least the understanding of
words \emph{essentially different} and \emph{function theory} should
be clarified. We discuss approaches to do that. We also present a new
framework for analytic function theories based on group 
representations.
\end{abstract}
\section{Introduction}
The classic heritage of complex analysis is contested between several
complex variables theory and hypercomplex analysis. The first one
was founded long ago by Cauchy and Weierstrass themselves and sometime
thought to be the only crown-prince. The hypercomplex analysis is not
a single theory but a family of related constructions discovered quite
recently~\cite{BraDelSom82,DelSomSou92,GuerlebeckSproessig90} 
(and rediscovered up to now)
under hypercomplex framework.

Such a variety of theories puts the question on their classification.
One could dream about a Mendeleev-like periodic table for hypercomplex
analysis, which clearly explains properties of different theories,
relationship between them and indicates how many blank cells are
waiting for us. Moreover, because hypercomplex analysis is the
recognized background for classic and quantum theories like the
Maxwell and Dirac equations, such a
table could play the role of \emph{the Mendeleev table for elementary
particles and fields}. We will return to this metaphor
and find it is not very superficial.

To make a step in the desired direction we should specify the notion
of \emph{function theory} and define the concept of \emph{essential
difference}. Probably many people agree that the core of complex
analysis consist of
\begin{enumerate}
\item\label{it:first} The Cauchy-Riemann equation and complex
derivative $\frac{\partial }{\partial z}$;
\item The Cauchy theorem;
\item The Cauchy integral formula;
\item The Plemeli-Sokhotski formula;
\item\label{it:last} The Taylor and Laurent series.
\end{enumerate}
Any development of several complex variables theory or hypercomplex
analysis is beginning from analogies to these notions and results. 
Thus we adopt the following
\begin{defn}
A \emph{function theory} is a collection of notions and results, which
includes at least analogies of~\ref{it:first}--\ref{it:last}.
\end{defn}
Of course the definition is more philosophical than
mathematical. For example, the understanding of an \emph{analogy} and
especially the \emph{right} analogy usually generates many disputes.

Again as a first approximation we propose the following
\begin{defn}
Two function theories is said to be \emph{similar} if there is a
correspondence between their objects such that analogies
of~\ref{it:first}--\ref{it:last} in one theory follow from their
counterparts in another theory. Two function theories are
\emph{essentially different} if they are not similar.
\end{defn}
Unspecified ``correspondence'' should probably be a linear map or
something else and we will look for its meaning soon. But it is
already clear that the \emph{similarity} is an equivalence relation
and we are looking for quotient sets with respect to it.

The layout is following. In Subsection~\ref{ss:factorization} the 
classic scheme of hypercomplex analysis is discussed and a possible 
variety of function theories appears. But we will see in 
Subsection~\ref{ss:reduce} that not all of them are very different.
Connection between group representations and (hyper)complex analysis 
is presented in Section~\ref{se:towards}. 
It could be a base for classification of essentially different 
theories.

\section{Abstract Nonsense about Function Theory}
In this Section we repeat shortly the scheme of development of
Clifford analysis as it could be found
in~\cite{BraDelSom82,DelSomSou92}. We examine different options
arising on this way and demonstrate that some differences are only
apparent not essential.

\subsection{Factorizations of the Laplacian}\label{ss:factorization}
We would like to see how the contents of~\ref{it:first}--\ref{it:last}
could be realized in a function theory. We are interested in function 
theories defined in $\Space{R}{d}$. The Cauchy theorem and
integral formula clearly indicates that the behavior of functions
inside a domain should be governed by their values on the boundary.
Such a property is particularly possessed by solutions to the
\emph{second order elliptic differential operator} $P$
\begin{displaymath}
P(x,\partial_x)=\sum_{i,j=1}^{d} a_{ij}(x)\partial_i \partial_j +
\sum_{i=1}^{d} b_i(x) \partial_i +c(x)
\end{displaymath}
with some special properties. Of course, the principal example is the
Laplacian
\begin{equation}\label{eq:Laplace}
\Delta=\sum_{i=1}^{d} \frac{\partial ^2}{\partial x_i^2}.
\end{equation}
\begin{enumerate}
\item \label{it:equation}\emph{Choice of different operators} (for
example, the Laplacian or the Helmholtz operator) is the first option which brings the
variety in the family of hypercomplex analysis.
\end{enumerate}

The next step is called \emph{linearization}. Namely we are looking
for two \emph{first order} differential operators $D$ and $D'$ such
that
\begin{displaymath}
DD'=P(x,\partial_x).
\end{displaymath}
The Dirac motivation to do that is to ``look for an equation linear in
in time derivative $\frac{\partial }{\partial t}$, because the
Schr\"odinger
equation is''. From the function theory point of view the 
Cauchy-Riemann operator should be linear also. But the most important 
gain of the step is an introduction of the Clifford algebra. For 
example, to factorize the Laplacian~\eqref{eq:Laplace} we put
\begin{equation}\label{eq:dirac}
D=\sum_{i=1}^d e_i \partial_i
\end{equation}
where $e_i$ are the \emph{Clifford algebra} generators:
\begin{equation}\label{eq:anti-comm}
e_i e_j + e_j e_i = 2\delta _{ij}, \qquad 1\leq i,j\leq d.
\end{equation}
\begin{enumerate}\addtocounter{enumi}{1}
\item \label{it:algebra}\emph{Different linearizations of a second
order operator} multiply the spectrum of theories.
\end{enumerate}
Mathematicians and physicists are looking up to now new factorization
even for the Laplacian\nocite{Keller93,ShaVas94a}. The essential
uniqueness of such factorization was already felt by Dirac himself but
it was never put as a theorem. So the idea of the \emph{genuine}
factorization becomes the philosophers' stone of our times.

After one made a choice~\ref{it:equation} and~\ref{it:algebra} the
following turns to be a routine. The equation
\begin{displaymath}
D' f(x)=0,
\end{displaymath}
plays the role of the Cauchy-Riemann equation.
Having a fundamental solution $F(x)$
to the operator $P(x,\partial _x)$ the Cauchy integral kernel defined
by
\begin{displaymath}
E(x)=D' F(x)
\end{displaymath}
with the property $DE(x)=\delta (x)$. Then the Stocks theorem implies
the Cauchy theorem and Cauchy integral formula. A decomposition of the
Cauchy kernel of the form
\begin{displaymath}
C(x-y)=\sum_{\alpha } V_\alpha (x) W_\alpha (y),
\end{displaymath}
where $V_\alpha (x)$ are some polynomials, yields via integration over
the ball the Taylor and Laurent series\footnote{Not all such 
decompositions give interesting series. The scheme from 
Section~\ref{se:towards} gives a selection rule to distinguish them.}. 
In such a way the program-minimum~\ref{it:first}--\ref{it:last} could 
be accomplished.

Thus all possibilities to alter function theory concentrated
in~\ref{it:equation} and \ref{it:algebra}. Possible universal algebras
arising from such an approach were investigated by
\person{F.~Sommen}~\cite{Sommen95a}. In spite of the apparent wide
selection, for operator $D$ and $D'$ with constant coefficients it was
found ``nothing dramatically new''~\cite{Sommen95a}:
\begin{quotation}
Of course one can study all these algebras and prove theorems or work
out lots of examples and representations of universal algebras. But in
the constant coefficient case the most important factorization seems
to remain the relation $\Delta=\sum x_j^2$, i.e., the one leading to
the definition of the Clifford algebra.
\end{quotation}
We present an example that there is no dramatical news not only on the
level of universal algebras but also for function theory (for the 
constant coefficient case). We will return to non
constant case in Section~\ref{se:towards}.

\subsection{Example of Connection}\label{ss:reduce}
We give a short example of similar theories with explicit
connection between them. The full account could be found
in~\cite{Kisil95c},
another example was considered in~\cite{Ryan95}.

Due to physical application we will consider equation
\begin{equation}\label{eq:mass}
\frac{\partial f}{\partial y_0}=(\sum_{j=1}^n e_j\frac{\partial
}{\partial y_j}+M)f,
\end{equation}
where $e_j$ are generators~\eqref{eq:anti-comm} of the Clifford 
algebra and $M=M_\lambda $ is an operator of
multiplication from the {\em right\/}-hand side by the Clifford number
$\lambda $. Equation~\eqref{eq:mass} is known in quantum
mechanics as the
{\em Dirac
equation for a particle with a non-zero rest
mass\/}~\cite[\S 20]{BerLif82}, \cite[\S 6.3]{BogShir80} and
\cite{Kravchenko95a}. We will specialize our results for the case
$M=M_\lambda$, especially for the simplest (but still important!) case
$\lambda\in\Space{R}{}$.
\begin{thm}
The function $f(y)$ is a solution to the equation
\begin{displaymath}
\frac{\partial f}{\partial y_0}=(\sum_{j=1}^n e_j\frac{\partial
}{\partial y_j}+M_1)f
\end{displaymath}
 if and only if the function
\begin{displaymath}
g(y)=e^{y_0 M_2} e^{-y_0 M_1} f(y)
\end{displaymath}
is a solution to the equation
\begin{displaymath}
\frac{\partial g}{\partial y_0}=(\sum_{j=1}^n e_j\frac{\partial
}{\partial y_j}+M_2)g,
\end{displaymath}
where $M_1$ and $M_2$ are bounded operators commuting with $e_j$.
\end{thm}
\begin{cor}\label{co:mass}
The function $f(y)$ is a solution to the equation~\eqref{eq:mass}
 if and only if the function
$e^{y_0 M}f(y)$
is a solution to the generalized Cauchy-Riemann
equation~\eqref{eq:dirac}.

In the case $M=M_\lambda$ we have $e^{y_0 M_\lambda}f(y)=f(y)e^{y_0
\lambda}$ and if $\lambda\in \Space{R}{}$ then $e^{y_0
M_\lambda}f(y)=f(y)e^{y_0 \lambda}=e^{y_0 \lambda}f(y)$.
\end{cor}

In this Subsection we construct a function theory (in the sense
of~\ref{it:first}--\ref{it:last}) for $M$-solutions of
the generalized Cauchy-Riemann operator based on Clifford analysis
and Corollary~\ref{co:mass}.

The set of solutions to~\eqref{eq:dirac} and~\eqref{eq:mass} in a nice
domain $\Omega$ will be denoted by
$\algebra{M}(\Omega)=\algebra{M}_0(\Omega)$ and
$\algebra{M}_M(\Omega)$ correspondingly. In the case $M=M_\lambda$ we
use the notation
$\algebra{M}_\lambda(\Omega)=\algebra{M}_{M_\lambda}(\Omega)$ also. We
suppose that all functions from $\algebra{M}_\lambda(\Omega)$ are
continuous in the closure of $\Omega$. Let
\begin{equation}\label{eq:cauchy-ker}
E(y-x)=
\frac{\Gamma(\frac{n+1}{2})}{2\pi^{(n+1)/2}}\,
\frac{\overline{y-x}}{\modulus{y-x}^{n+1}}
\end{equation}
be the Cauchy kernel~\cite[p.~146]{DelSomSou92}  and
\begin{displaymath}
d\sigma=\sum_{j=0}^n (-1)^j e_j dx_0 \wedge \ldots \wedge[dx_j] \wedge
\ldots \wedge dx_m.
\end{displaymath}
be the differential form of the ``oriented surface
element''~\cite[p.~144]{DelSomSou92}. Then for any
$f(x)\in\algebra{M}(\Omega)$ we have the Cauchy integral
formula~\cite[p.~147]{DelSomSou92}
\begin{displaymath}
\int_{\partial \Omega} E(y-x)\,d\sigma_y\,
f(y)=\left\{\begin{array}{cl}
f(x),& x\in\Omega\\
0,& x\not\in\bar{\Omega}
\end{array}.\right.
\end{displaymath}

\begin{thm}[Cauchy's Theorem]
Let $f(y)\in \algebra{M}_M(\Omega)$. Then
\begin{displaymath}
\int_{\partial \Omega} d\sigma_y\,e^{-y_0 M}f(y)=0.
\end{displaymath}
Particularly, for $f(y)\in \algebra{M}_\lambda(\Omega)$ we have
\begin{displaymath}
\int_{\partial \Omega} d\sigma_y\,f(y)e^{-y_0 \lambda}=0,
\end{displaymath}
and
\begin{displaymath}
\int_{\partial \Omega} d\sigma_ye^{-y_0 \lambda}\,f(y)=0,
\end{displaymath}
if $\lambda\in\Space{R}{} $.
\end{thm}

\begin{thm}[Cauchy's Integral Formula]
Let $f(y)\in \algebra{M}_M(\Omega)$. Then
\begin{equation}\label{eq:m-cauchy}
e^{x_0 M}\int_{\partial \Omega} E(y-x)\,d\sigma_y\,
e^{-y_0 M}f(y)=\left\{\begin{array}{cl}
f(x),& x\in\Omega\\
0,& x\not\in\bar{\Omega}
\end{array}.\right.
\end{equation}
Particularly, for $f(y)\in \algebra{M}_\lambda(\Omega)$ we have
\begin{displaymath}
\int_{\partial \Omega} E(y-x)\,d\sigma_y\,
f(y)e^{(x_0-y_0) \lambda}=\left\{\begin{array}{cl}
f(x),& x\in\Omega\\
0, &x\not\in\bar{\Omega}
\end{array}.\right.
\end{displaymath}
and
\begin{displaymath}
\int_{\partial \Omega} E(y-x)e^{(x_0-y_0) \lambda}\,d\sigma_y\,
f(y)=\left\{\begin{array}{cl}
f(x),& x\in\Omega\\
0,& x\not\in\bar{\Omega}
\end{array}.\right.
\end{displaymath}
if $\lambda\in\Space{R}{} $.
\end{thm}
It is hard to expect that formula~\eqref{eq:m-cauchy} may be rewritten
as
\begin{displaymath}
\int_{\partial \Omega} E'(y-x)\,d\sigma_y\,
f(y)=\left\{\begin{array}{cl}
f(x),& x\in\Omega\\
0,& x\not\in\bar{\Omega}
\end{array}\right.
\end{displaymath}
with a simple function $E'(y-x)$.

Because an application of the bounded operator $e^{y_0 M}$ does not
destroy uniform convergency of functions we obtain
(cf.~\cite[Chap.~II, \S~0.2.2, Theorem~2]{DelSomSou92})
\begin{thm}[Weierstrass' Theorem]
Let $\{f_k\}_{k\in\Space{N}{}}$ be a sequence in
$\algebra{M}_M(\Omega)$, which converges uniformly to $f$ on each
compact subset $K\in \Omega$. Then
\begin{enumerate}
\item $f\in \algebra{M}_M(\Omega)$.
\item For each multi-index
$\beta=(\beta_0,\ldots,\beta_m)\in\Space{N}{n+1}$, the sequence
$\{\partial ^\beta f_k\}_{k\in\Space{N}{}} $ converges uniformly on
each compact subset $K\in \Omega$ to $\partial ^\beta f$.
\end{enumerate}
\end{thm}

\begin{thm}[Mean Value Theorem]
Let $f\in\algebra{M}_M(\Omega)$. Then for all $x\in \Omega$ and $R>0$
such that the ball $\Space{B}{}(x,R)\in\Omega $,
\begin{displaymath}
f(x)= e^{x_0 M}
\frac{(n+1)\Gamma(\frac{n+1}{2})}{2R^{n+1}\pi^{(n+1)/2}}
\int_{\Space{B}{}(x,R)} e^{-y_0 M} f(y)\,dy.
\end{displaymath}
\end{thm}

Such a reduction of theories could be pushed even 
future~\cite{Kisil95c} up to the notion of 
hypercomplex differentiability~\cite{Malonek93}, but we will 
stop here.  

\section{Hypercomplex Analysis and Group
Representations --- Towards a Classification}\label{se:towards}
To construct a classification of non-equivalent objects one could
use their groups of symmetries. Classical example is Poincar\'e's proof
of bi-holomorphic non-equivalence of the unit ball and polydisk via
comparison their groups of bi-holomorphic automorphisms. To employ 
this approach we need a construction of hypercomplex analysis from its
symmetry group. The following scheme is firstly presented here (up to 
the author knowledge) and has its roots 
in~\cite{Kisil94e,Kisil95d,Kisil95a}.

Let $G$ be a group which acts via transformation of a closed domain
$\bar{\Omega}$. Moreover, let $G: \partial \Omega\rightarrow \partial
\Omega$ and $G$ act on $\Omega$ and $\partial \Omega$ transitively.
Let us fix a point $x_0\in \Omega$ and let $H\subset G$ be a
stationary subgroup of point $x_0$. Then domain $\Omega$ is naturally
identified with the  homogeneous space $G/H$. Till the moment we do
not request anything untypical. Now let 
\begin{itemize}
\item\emph{there exist a $H$-invariant measure $d\mu$ on $\partial
\Omega$}.
\end{itemize}
 We consider the Hilbert space $\FSpace{L}{2}(\partial
\Omega, d\mu)$. Then geometrical transformations of $\partial \Omega$
give us the representation $\pi$ of $G$ in $\FSpace{L}{2}(\partial
\Omega, d\mu)$.
 Let $f_0(x)\equiv 1$ and $\FSpace{F}{2}(\partial
\Omega, d\mu)$ be the closed liner subspace of $\FSpace{L}{2}(\partial
\Omega, d\mu)$ with the properties:
\begin{enumerate}
\item\label{it:begin} $f_0\in \FSpace{F}{2}(\partial \Omega, d\mu)$;
\item $\FSpace{F}{2}(\partial \Omega, d\mu)$ is $G$-invariant;
\item\label{it:end} $\FSpace{F}{2}(\partial \Omega, d\mu)$ is $G$-irreducible.
\end{enumerate}
The \emph{standard wavelet transform} $W$ is defined by
\begin{displaymath}
W: \FSpace{F}{2}(\partial \Omega, d\mu) \rightarrow
\FSpace{L}{2}(G): f(x) \mapsto
\widehat{f}(g)=\scalar{f(x)}{\pi(g)f_0(x)}_{\FSpace{L}{2}(\partial
\Omega,d\mu) }
\end{displaymath}
Due to the property $[\pi(h)f_0](x)=f_0(x)$, $h\in H$ and 
identification $\Omega\sim G/H$ it could be translated to the embedding:
\begin{equation}\label{eq:cauchy}
\widetilde{W}: \FSpace{F}{2}(\partial \Omega, d\mu) \rightarrow
\FSpace{L}{2}(\Omega): f(x) \mapsto
\widehat{f}(y)=\scalar{f(x)}{\pi(g)f_0(x)}_{\FSpace{L}{2}(\partial
\Omega,d\mu) },  
\end{equation}
where $y\in\Omega$ for some $ h\in H$. The imbedding~\eqref{eq:cauchy} 
is \emph{an abstract analog of the Cauchy integral formula}. 
Let functions $V_\alpha $ be the \emph{special functions} generated by 
the representation of $H$. Then the decomposition of 
$\widehat{f}_0(y)$ by $V_\alpha $ gives us the Taylor series.

The scheme is inspired by the following interpretation of complex 
analysis.
\begin{example}
Let the domain $\Omega$ be the unit disk $\Space{D}{}$, $\partial 
\Space{D}{}=\Space{S}{}$. We select
the group $SL(2,\Space{R}{})\sim SU(1,1)$ acting on $\Space{D}{} $
via the fractional-linear transformation:
\begin{displaymath}
\matr{a}{b}{c}{d}: z\mapsto \frac{az+b}{cz+d}.
\end{displaymath}
We fix $x_0=0$. Then its stationary group is $U(1)$ of rotations of 
$\Space{D}{}$. Then the Lebesgue measure on $\Space{S}{} $ is 
$U(1)$-invariant. We obtain $\Space{D}{}\sim SL(2,\Space{R}{} )/U(1)$. 
The subspace of $\FSpace{L}{2}(\Space{S}{},dt)$ satisfying 
to~\ref{it:begin}--\ref{it:end} is the Hardy space. The wavelets 
transform\eqref{eq:cauchy} give exactly the Cauchy formula.
The proper functions of $U(1)$ are exactly $z^n$, which provide the 
base for the Taylor series. The Riemann mapping theorem allows to 
apply the scheme to any connected, simply-connected domain.
\end{example}

The conformal group of the M\"obius transformations plays the same 
role in Clifford analysis. One usually says that the conformal group 
in $\Space{R}{n}$, $n>2$ is
not so rich as the conformal group in $\Space{R}{2}$.
Nevertheless, the conformal covariance has many applications in
Clifford analysis~\cite{Cnops94a,Ryan95b}.
Notably,
groups of conformal mappings of open unit balls $\Space{B}{n} \subset
\Space{R}{n}$ \emph{on}to itself are similar for all $n$ and
as sets can be
parametrized by the product of \Space{B}{n} itself and the group
of isometries of its boundary \Space{S}{n-1}.
\begin{thm}{\cite{Kisil95i} }\label{th:ball}
Let $a\in\Space{B}{n}$, $b\in\Gamma_n$ then the M\"obius
transformations of
the form
\begin{displaymath}
\phi_{(a,b)}=\matr{b}{0}{0}{{b}^{*-1}}\matr{1}{-
a}{{a}^*}{-1}=\matr{b}{-ba}{{b}^{*-1}a^*}
{-{b}^{*-1}},
\end{displaymath}
constitute the group $B_{n}$ of conformal mappings of the open unit
ball $\Space{B}{n}$ onto itself. $B_{n}$ acts on
$\Space{B}{n}$ transitively.
Transformations of the form $\phi_{(0,b)}$  constitute a
subgroup isomorphic to $\object[(n)]{O}$. The homogeneous space
$B_{n}/\object[(n)]{O}$ is isomorphic as a set to
$\Space{B}{n}$. Moreover:
\begin{enumerate}
\item $\phi_{(a,1)}^2=1$ identically on $\Space{B}{n}$
($\phi_{(a,1)}^{-1}=\phi_{(a,1)}$).
\item $\phi_{(a,1)}(0)=a$, $\phi_{(a,1)}(a)=0$.
\end{enumerate}
\end{thm}

Obviously, conformal mappings preserve the space of null solutions to
the \emph{Laplace} operator~\eqref{eq:Laplace} and null solutions the
\emph{Dirac} operator~\eqref{eq:dirac}. The
group $B_{n}$ is sufficient for construction of the Poisson and the 
Cauchy integral representation of harmonic functions and Szeg\"o and 
Bergman projections in Clifford analysis by the 
formula~\cite{Kisil95d}
\begin{equation}\label{eq:reproduce}
K(x,y)=c\int_G
[\pi_g f](x) \overline{[\pi_g f](y)}\,dg,
\end{equation}
where $\pi_g$ is an irreducible unitary
square integrable representation of a group $G$,
$f(x)$ is an arbitrary non-zero function, and $c$ is a constant.

The scheme gives a correspondence between \emph{function theories} and 
\emph{group representations}. The last are rather well studded and 
thus such a connection could be a foundation for a classification of 
function theories. Particularly, the \emph{constant coefficient} 
function theories in the sense of \person{F.~Sommen}\cite{Sommen95a} 
corresponds to the groups acting only on the function domains in the 
Euclidean space. Between such groups the Moebius transformations play 
the leading role. On the contrary, the \emph{variable coefficient} 
case is described by groups acting on the function space in the 
non-point sense (for example, combining action on the functions domain 
and range, see~\cite{Kisil94e}). The set of groups of the second kind 
should be more profound.

\begin{rem}
It is known that many results in real analysis~\cite{McIntosh95a} several 
variables theory~\cite{MiSha95} could be obtained or even explained via 
hypercomplex analysis. One could see roots of this phenomenon in 
relationships between groups of geometric symmetries of two theories: the 
group of hypercomplex analysis is wider.
\end{rem}

Returning to our metaphor on the Mendeleev table we would like recall
that it began as linear ordering with respect to atomic masses but
have received an explanation only via representation theory of the
rotation group.

Our consideration provides a ground for the following
\begin{conj}
Most probably there is the only constant coefficient function theory
on the Euclidean space or at most there are two of them.
\end{conj}

\section{Acknowledgments and Apologies}
It is my pleasure to express my gratitude to the J.~Cnops, 
R.~Delanghe, K.~Guerlebeck, V.V.~Kravchenko, I.~Mitelman, J.~Ryan, 
M.~Shapiro, F.~Sommen, W.~Spr\"o\-{\ss}ig, N.~Vasilevski who share 
with me their ideas on hypercomplex analysis. It was especially 
stimulating for this paper author's stay at Universiteit Gent, 
Vakgroep Wiskundige Analyse (Belgium) under INTAS grant 93--0322. I am 
also grateful to DFG for the financial support of my participation at 
the Conference.

The bibliography to such a paper should be definitely more complete 
and representative. Unfortunately, I mentioned only a few papers 
among deserving it.

%\bibliographystyle{abbrv}
%\bibliography{MRABBREV,akisil,analyse,aclifford,aphysics}

\begin{thebibliography}{10}

\bibitem{BerLif82}
V.~B. Berestetskii, E.~M. Lifshitz, and L.~P. Pitaevski.
\newblock {\em Quantum Electrodynamics}, volume~4 of {\em L. D.~Landay and E.
  M.~Lifshitz, Course of Theoretical Physics}.
\newblock Pergamon Press, Oxford, second edition, {\noopsort{}}1982.

\bibitem{BogShir80}
N.~N. Bogoliubov and D.~V. Shirkov.
\newblock {\em Introduction to the Theory of Quantized Fields}.
\newblock John Willey \& Sons, Inc., New York, third edition,
  {\noopsort{}}1980.

\bibitem{BraDelSom82}
F.~Brackx, R.~Delanghe, and F.~Sommen.
\newblock {\em Clifford Analysis}, volume~76 of {\em Research Notes in
  Mathematics}.
\newblock Pitman (Advanced Publishing Program), Boston, MA, 1982.

\bibitem{Cnops94a}
J.~Cnops.
\newblock {\em {Hurwitz} Pairs and Applications of {M\"obius} Transformations}.
\newblock {Habilitation} dissertation, Faculteit van de Wetenschappen,
  Universiteit Gent, 1994.
\newblock See also~\cite{Cnops02a}.

\bibitem{Cnops02a}
J.~Cnops.
\newblock {\em An Introduction to {D}irac Operators on Manifolds}, volume~24 of
  {\em Progress in Mathematical Physics}.
\newblock Birkh\"auser Boston Inc., Boston, MA, 2002.

\bibitem{DelSomSou92}
R.~Delanghe, F.~Sommen, and V.~Sou{\v{c}}ek.
\newblock {\em Clifford Algebra and Spinor-Valued Functions. {A} function
  theory for the {Dirac} operator}, volume~53 of {\em Mathematics and its
  Applications}.
\newblock Kluwer Academic Publishers Group, Dordrecht, 1992.
\newblock Related REDUCE software by F. Brackx and D. Constales, With 1 IBM-PC
  floppy disk (3.5 inch).

\bibitem{GuerlebeckSproessig90}
K.~G{\"u}rlebeck and W.~Spr{\"o}ssig.
\newblock {\em Quaternionic Analysis and Elliptic Boundary Value Problems},
  volume~89 of {\em International Series of Numerical Mathematics}.
\newblock Birkh\"auser Verlag, Basel, 1990.

\bibitem{Keller93}
J.~Keller.
\newblock The geometric content of the electron theory. part {I}.
\newblock {\em Adv. in Appl. Clifford Algebras}, 3(2):147--200, 1993.

\bibitem{Kisil95c}
V.~V. Kisil.
\newblock Connection between different function theories in {C}lifford
  analysis.
\newblock {\em Adv. Appl. Clifford Algebras}, 5(1):63--74, 1995.
\newblock \arXiv{funct-an/9501002}.

\bibitem{Kisil95a}
V.~V. Kisil.
\newblock Integral representations and coherent states.
\newblock {\em Bull. Belg. Math. Soc. Simon Stevin}, 2(5):529--540, 1995.
\newblock
  \href{http://projecteuclid.org/getRecord?id=euclid.bbms/1103408676}{On-line}.

\bibitem{Kisil95d}
V.~V. Kisil.
\newblock Construction of integral representations for spaces of analytic
  functions.
\newblock {\em Dokl. Akad. Nauk}, 350(4):446--448, 1996.
\newblock (Russian) \MR{98d:46027}.

\bibitem{Kisil95i}
V.~V. Kisil.
\newblock M\"obius transformations and monogenic functional calculus.
\newblock {\em Electron. Res. Announc. Amer. Math. Soc.}, 2(1):26--33, 1996.
\newblock
  \href{http://www.ams.org/era/1996-02-01/S1079-6762-96-00004-2/}{On-line}.

\bibitem{Kisil94e}
V.~V. Kisil.
\newblock Relative convolutions. {I}. {P}roperties and applications.
\newblock {\em Adv. Math.}, 147(1):35--73, 1999.
\newblock \arXiv{funct-an/9410001},
  \href{http://www.idealibrary.com/links/doi/10.1006/aima.1999.1833}{On-line}.
  \Zbl{933.43004}.

\bibitem{Kravchenko95a}
V.~V. Kravchenko.
\newblock On biquaternionic bag model.
\newblock {\em Zeitschrift f\"ur Anal. und ihre Anwend.}, 14(1):3--14, 1995.

\bibitem{Malonek93}
H.~R. Malonek.
\newblock Hypercomplex differentiability and its applications.
\newblock In F.~Brackx, R.~Delanghe, and H.~Serras, editors, {\em {C}lifford
  Algebras and Their Applications in Mathematical Physics, Proceedings of the
  {T}hird International Conference held in {D}einze}, volume~55 of {\em
  Fundamental Theories of Physics}, pages 141--150. Kluwer Academic Publishers
  Group, Dordrecht, 1993.
\newblock \MR{94j:00019}.

\bibitem{McIntosh95a}
A.~McIntosh.
\newblock {Clifford} algebras, {Fourier} theory, singular integral operators,
  and partial differential equations on {Lipschitz} domains.
\newblock In Ryan \cite{Cliff95}, pages 33--88.

\bibitem{MiSha95}
I.~M. Mitelman and M.~V. Shapiro.
\newblock Differentiation of the {M}artinelli-{B}ochner integrals and the
  notion of hyperderivability.
\newblock {\em Math. Nachr.}, 172:211--238, 1995.

\bibitem{Cliff95}
J.~Ryan, editor.
\newblock {\em {Clifford} Algebras in Analysis and Related Topics}.
\newblock CRC Press, Boca Raton, {\noopsort{}}1995.

\bibitem{Ryan95b}
J.~Ryan.
\newblock Some application of conformal covariance in {Clifford} analysis.
\newblock In Ryan95b \cite{Cliff95}, pages 128--155.

\bibitem{Ryan95}
J.~Ryan.
\newblock Intrinsic {Dirac} operator in{ \Space{C}{n}}.
\newblock {\em Adv. in Math.}, 118(1):99--133, 1996.

\bibitem{ShaVas94a}
M.~Shapiro and N.~Vasilevski.
\newblock Quaternionic $\psi$-hyperholomorphic functions, singular integral
  operators and boundary value problems. {I}. $\psi$-{H}yperholomorphic
  function theory.
\newblock {\em Complex Variables Theory Appl.}, 27, 1994.

\bibitem{Sommen95a}
F.~Sommen.
\newblock Universal algebras coming from factorization.
\newblock page~10, 1995.
\newblock Preprint.

\end{thebibliography}
\providecommand{\noopsort}[1]{} \providecommand{\printfirst}[2]{#1}
  \providecommand{\singleletter}[1]{#1} \providecommand{\switchargs}[2]{#2#1}
  \providecommand{\irm}{\textup{I}} \providecommand{\iirm}{\textup{II}}
  \providecommand{\vrm}{\textup{V}} \providecommand{\cprime}{'}
  \providecommand{\eprint}[2]{\texttt{#2}}
  \providecommand{\myeprint}[2]{\texttt{#2}}
  \providecommand{\arXiv}[1]{\myeprint{http://arXiv.org/abs/#1}{arXiv:#1}}
  \providecommand{\doi}[1]{\href{http://dx.doi.org/#1}{doi: #1}}

\end{document}
\bibitem{Yeh}
R.Z.~Yeh.
\newblock Analysis and applications of holomorphic functions 
in higher dimensions. 
\newblock {\em Transactions of the AMS\/}, 345(1):151--177, 1994.